\documentclass[11pt]{article}

\usepackage{amsmath}
\usepackage{amssymb}
\usepackage{eucal}

\begin{document}

\baselineskip 16pt

\title{On  generalized $\sigma$-soluble groups\thanks{Research is supported by an NNSF grant of China (Grant No.
11401264) and a TAPP of Jiangsu Higher Education Institutions (PPZY 2015A013)}}

\author{Jianhong Huang\\\
{\small School of Mathematics and Statistics, Jiangsu Normal University,}\\ {\small Xuzhou, 221116, P.R. China}\\
{\small E-mail: jhh320@126.com}\\ \\
Bin Hu    \thanks{Corresponding author}\\
{\small School of Mathematics and Statistics, Jiangsu Normal University,}\\
{\small Xuzhou 221116, P. R. China}\\
{\small E-mail: hubin118@126.com}\\ \\
{ Alexander  N. Skiba   }\\
{\small Department of Mathematics and Technologies of Programming,  Francisk Skorina Gomel State University,}\\
{\small Gomel 246019, Belarus}\\
{\small E-mail: alexander.skiba49@gmail.com}}

\date{}                               
\maketitle

\begin{abstract} Let $\sigma =\{\sigma_{i} | i\in I\}$ be a partition of
 the set of all primes $\Bbb{P}$ and  $G$  a finite group.   Let $\sigma 
(G)=\{\sigma _{i} : \sigma _{i}\cap   \pi (G)\ne \emptyset$. 
A set  ${\cal H}$  of subgroups of $G$ is said to be a  \emph{complete Hall
 $\sigma $-set} of $G$   if   every  member $\ne 1$ of
 ${\cal H}$ is a Hall $\sigma _{i}$-subgroup  of $G$ for some $i\in I$ and
 $\cal H$ contains exactly one Hall  $\sigma _{i}$-subgroup of $G$ for 
every  $i$ such that $\sigma _{i}\in \sigma (G)$. 
 We say that $G$ is  \emph{$\sigma$-full}  if $G$ possesses a  complete Hall
 $\sigma $-set. 
 A complete Hall $\sigma $-set $\cal H$ of $G$  is said to be a 
\emph{$\sigma$-basis}  of $G$ if every two subgroups $A, B \in\cal H$ are 
permutable, that is, $AB=BA$.

 In this paper, we study properties of finite groups having a $\sigma$-basis.
   In particular,  we prove that if $G$    has a  a $\sigma$-basis, then $G$ is
 \emph{generalized $\sigma$-soluble}, that is, $G$ has a complete Hall
 $\sigma $-set  and for every chief factor $H/K$ of 
$G$ we have $|\sigma (H/K)|\leq 2$. Moreover, answering to Problem 8.28 in \cite{commun},
we prove the following

{\bf Theorem A.} {\sl Suppose that  $G$ is $\sigma$-full. Then every complete
 Hall $\sigma$-set of $G$ forms a $\sigma$-basis of $G$ if and only if $G$ is
 generalized $\sigma$-soluble and  for the automorphism group $G/C_{G}(H/K)$,
  induced by $G$ on any  its 
chief factor $H/K$,  we have either   $\sigma (H/K)=\sigma (G/C_{G}(H/K))$ 
 or  $\sigma (H/K) =\{\sigma _{i}\}$  and  $G/C_{G}(H/K)$ is a 
$\sigma 
_{i} \cup   \sigma _{j}$-group  for some $i\ne j$. 
}

  \end{abstract}

\footnotetext{Keywords: finite group,  Hall subgroup, $\sigma$-semipermutable  
subgroup, $\sigma$-basis,  generalized ${\sigma}$-soluble   group.}

\footnotetext{Mathematics Subject Classification (2010): 20D10, 20D15}
\let\thefootnote\thefootnoteorig

\section{Introduction}

Throughout this paper, all groups are finite and $G$ always denotes
a finite group. Moreover,  $\mathbb{P}$ is the set of all  primes,
  $\pi= \{p_{1}, \ldots , p_{n}\} \subseteq  \Bbb{P}$ and  $\pi' =  \Bbb{P} \setminus \pi$.
 If
 $n$ is an integer, the symbol $\pi (n)$ denotes
 the set of all primes dividing $n$; as usual,  $\pi (G)=\pi (|G|)$, the set of all
  primes dividing the order of $G$.

In what follows, $\sigma$  is some partition of  
$\Bbb{P}$, that is,  $\sigma =\{\sigma_{i} |
 i\in I \}$, where   $\Bbb{P}=\bigcup_{i\in I} \sigma_{i}$
 and $\sigma_{i}\cap
\sigma_{j}= \emptyset  $ for all $i\ne j$.   
By the analogy with the notation   $\pi (n)$, we write  $\sigma (n)$ to denote 
the set  $\{\sigma_{i} |\sigma_{i}\cap \pi (n)\ne 
 \emptyset  \}$;   $\sigma (G)=\sigma (|G|)$.

 A \emph{complete  set of Sylow subgroups $\cal S$ of $G$}  
contains exactly one Sylow $p$-subgroup for each prime $p$, that is,
 ${\cal S}=\{1, P_{1}, \ldots P_{t}\}$, where $P_{i}$ is a Sylow $p_{i}$-subgroup of $G$ and 
$\pi (G)=\{p_{1}, \ldots , p_{t}\}$.
The set  $\cal S$  is said to be a \emph{ Sylow basis of $G$} provided every  two 
subgroups  $P, Q \in {\cal S}$ are permutable, that is, $PQ=QP$.

In general, we say that a   set  ${\cal H}$ of subgroups of $G$ is a
 \emph{complete Hall $\sigma $-set} of $G$ \cite{commun}  if
 every member $\ne 1$ of  ${\cal H}$ is a Hall $\sigma _{i}$-subgroup of $G$
 for some $\sigma _{i} \in \sigma$ and ${\cal H}$ contains exactly one Hall
 $\sigma _{i}$-subgroup of $G$ for every  $\sigma _{i}\in  \sigma (G)$.  
We say that  $G$ is   \emph{$\sigma$-full} \cite{commun} if $G$ possesses  a complete Hall 
$\sigma $-set.  A complete Hall $\sigma $-set $\cal H$ of $G$  is said to be a 
\emph{$\sigma$-basis}  of $G$ \cite{2} if every two subgroups $A, B \in\cal H$ are 
permutable, that is, $AB=BA$.

In this paper we deal with the following  two generalizations of 
solubility.

{\bf Definition 1.1.} We say that  $G$ is:

(i)  \emph{$\sigma$-soluble} 
\cite{2} if  $|\sigma (H/K)|= 1$ for every  chief factor $H/K$ of $G$.

(i)  \emph{generalized$\sigma$-soluble} 
 if  $G$ is  $\sigma$-full and   $|\sigma (H/K)|\leq  2$ for every
  chief factor $H/K$ of $G$.

 Before continuing, consider three classical  cases. 

{\bf Example 1.2.} (i) In the classical case,
 when $\sigma =\sigma ^{0}=\{\{2\}, \{3\}, \ldots 
\}$,  $G$ is     $\sigma ^{0}$-soluble  if
 and only if it is soluble.  Note also that in view of the Burnside's 
$p^{a}q^{b}$-theorem, $G$ is   $\sigma ^{0}$-soluble  if
 and only if it  is  generalized $\sigma ^{0}$-soluble.

(ii) In the other classical case, when  $\sigma =\sigma ^{\pi}=\{\pi, 
\pi'\}$  $G$  is: $\sigma ^{\pi}$-soluble  if   and only if $G$ is 
   $\pi$-separable;  generalized $\sigma ^{\pi}$-soluble  if   and only if 
$G$ has both  a Hall $\pi$-subgroup  and a  Hall $\pi'$-subgroup.

(iii) In fact, in the theory of $\pi$-soluble groups ($\pi= \{p_{1}, \ldots , p_{n}\}$)
 we deal with the  partition 
$\sigma =\sigma ^{0\pi }=\{\{p_{1}\}, \ldots , \{p_{n}\}, \pi'\}$ of $\Bbb{P}$. 
Note that  $G$ is $\sigma ^{0\pi }$-soluble  if and only if it is 
   $\pi$-soluble.  Note also that $G$ is generalized   $\sigma ^{0\pi 
}$-soluble  if and only if the following hold: (a) $G$ has both
  a Hall $\pi$-subgroup $E$  and a  Hall $\pi'$-subgroup, and  (b) $(H\cap E)K/K$ is a
 Sylow subgroup of $H/K$   for   every chief factor $H/K$ of $G$.

The well-known  Hall's theorem  states that  $G$ is
 soluble if and only if it has a Sylow basis.  This classical result makes 
natural to ask: {Does it true that $G$ is
 $\sigma$-soluble if and only if it has a $\sigma$-basis?}

A partial confirmation of this hypothesis is given by the following

{\bf Theorem 1.3} (See Theorem A in \cite{2}).  {\sl Every $\sigma$-soluble group 
possesses a  $\sigma$-basis}

Now consider the following 

{\bf Example 1.4.}   Let $\sigma =\{\{2, 3\}, \{5\}, \{2, 3, 5\}'\}$, and
 let $G=A_{5}\times (C_{11}\rtimes \text{Aut} (C_{11}))$, where $A_{5}$ is the 
alternating group of degree 5 and $C_{11}$ is a group of order 11.  Let 
$A_{4}\simeq A\leq A_{5}$ and $B$ a Sylow 5-subgroup of $A_{5}$. 
Let $H_{1}=AC_{2}$, $H_{2}=BC_{5}$ and $H_{3}=C_{11}$,
 where $C_{5}\times C_{2}=\text{Aut} 
(C_{11}))$.   
 Then   the set
 $\{H_{1}, H_{2}, H_{3}\}$ is  a $\sigma$-basis of $G$.  Nevertheless,  
$G$ is not   $\sigma$-soluble.  This example shows that in general the 
answer to above question is negative. 

 Nevertheless, the following result 
is true.

{\bf Theorem A.} {\sl Suppose that  $G$ possesses  a complete Hall $\sigma 
$-set $\cal H$.  }

(i)  {\sl If $\cal H$ is  a  $\sigma$-basis
 of $G$, then $G$  is generalized  $\sigma$-soluble. }

(ii)  {\sl If  $H\leq A\in {\cal H}$ and  
$HV^{x}=V^{x}H$ for all $x\in G$ and all $V\in \cal H$ such that
$(|H|, |V|)=1$, then $H^{G}$ is generalized $\sigma$-soluble. 
 
}

In view of the Burnside's $p^{a}q^{b}$-theorem, in the case where 
 $\sigma =\{\{2\}, \{3 \}, \ldots \}$ we get from Theorem A(ii) the following

{\bf Corollary 1.5} (Isaaks \cite{isaaks}). {\sl  If a $p$-subgroup $H$
 of $G$ permutes with all Sylow subgroups $P$ of $G$ such that $(p, |P|)=1$, then $H^{G}$ is soluble.   }

{\bf Corollary 1.4.} {\sl  If  $H\in {\cal H}$ and $HV^{x}=V^{x}H$ for all $x\in G$ and all $V\in \cal H$ such that
$(|H|, |V|)=1$, then $G$ is $\sigma _{i}'$-soluble where
 $\sigma (H)= \{\sigma _{i}\}$. 
}

{\bf Corollary 1.6} (Borovikov \cite{borov}). {\sl  If a Sylow $p$-subgroup $P$
 of $G$  is permutes with all Sylow subgroups $Q$ of $G$ such that $(p, |Q|)=1$, then $G$
 is $p'$-soluble.   }

The integers $n$ and $m$ are called \emph{$\sigma$-coprime}  if $\sigma (n)
 \cap \sigma (m)=\emptyset$.

In fact,  Theorem A(i) is a corollary of the following fact.

{\bf Theorem B.}    {\sl   Suppose that   $G=A_{1}A_{2}=A_{2}A_{3}=A_{1}A_{3},$  
where  $A_{1}$ is  $\sigma$-soluble and  $A_{2}$ and $A_{3}$
 are generalized $\sigma$-soluble subgroups of $G$.
 If for some $i, j, k\in I$  the  three
 indices
$|G:N_{G}(O^{\sigma _{i}}(A_{1}))|$, $|G:N_{G}(O^{\sigma _{j}}(A_{2}))|$,
 $|G:N_{G}(O^{\sigma _{k}}(A_{3}))|$ are pairwise  $\sigma$-coprime,
  then
$G$ is generalized $\sigma$-soluble. }

{\bf Corollary  1.7.}    {\sl   Suppose that   $A_{1}$ is a  
  $\sigma$-soluble subgroup  and  $A_{2}$ and $A_{3}$
 are generalized $\sigma$-soluble subgroups of $G$.
 If   the  three
 indices
$|G:A_{1}|$, $|G:A{2}|$,
 $|G:A_{3}|$ are pairwise  $\sigma$-coprime,
  then
$G$ is generalized $\sigma$-soluble. }

{\bf Corollary  1.8.}    {\sl  
 Suppose that $G=A_{1}A_{2}=A_{2}A_{3}=A_{1}A_{3},$  
where  $A_{1}$, $A_{2}$ and $A_{3}$ are soluble subgroups of $G$.
 If  the  three indices
$|G:N_{G}(A_{1})|$, $|G:N_{G}(A_{2})|$, $|G:N_{G}(A_{3})|$ are pairwise  coprime, then
$G$ is soluble. }

{\bf Corollary 1.8.} (H. Wielandt).  {\sl  If $G$ has
three soluble subgroups $A_{1}$, $A_{2}$ and $A_{3}$ whose indices
$|G:A_{1}|$, $|G:A_{2}|$, $|G:A_{3}|$ are pairwise  coprime, then
$G$ is itself soluble. }

Following \cite{commun}, we use $\frak{H}_{\sigma}$
to denote the class of all $\sigma$-full groups  $G$ such that every
complete Hall $\sigma$-set of $G$ forms a $\sigma$-basis of $G$.

In  \cite[VI, Section 3]{hupp}, Huppert described soluble groups $G$ in which 
every complete Sylow set of $G$  forms a Sylow basis of $G$. The results 
in \cite[VI, Section 3]{hupp} are motivations for 
    the following two questions.

{\bf Question 1.9} (See Problem 8.29 in \cite{commun}).  {\sl  Describe groups in
  $\frak{H}_{\sigma}$.}

{\bf  Question 1.10}   (See Problem 8.30 in \cite{commun}).    {\sl  Describe
 $\sigma$-soluble groups groups in  $\frak{H}_{\sigma}\cap  \frak{S}_{\sigma}$
}

As one more  application of Theorem A, we prove the
 following result, which gives the answer to Question 1.9.

{\bf Theorem C.} {\sl Suppose that  $G$ is $\sigma$-full. Then every complete
 Hall $\sigma$-set of $G$ forms a $\sigma$-basis of $G$ if and only if $G$ is
 generalized $\sigma$-soluble and  for the automorphism group $G/C_{G}(H/K)$,
  induced by $G$ on any  its 
chief factor $H/K$,   we have $|\sigma (G/C_{G}(H/K))|\leq 2$ and also 
 $\sigma(H/K)\subseteq  \sigma (G/C_{G}(H/K))$ in the case $|\sigma (G/C_{G}(H/K))|= 2$.
}

ON the base of Theorem C  we prove the following result which gives the 
answer to Question 1.10.

{\bf Theorem D.} {\sl  Suppose that $G$ is $\sigma$-soluble.
 Then the following statements are equivalent:}

(i) {\sl  Every complete Hall $\sigma$-set of $G$ forms a $\sigma$-basis of $G$. }

(ii) {\sl  The automorphism group $G/C_{G}(H/K)$,  induced by $G$ on any its chief factor
$H/K$ with $\sigma (H/K)=\{\sigma _{i}\}$,   is either  a  $\sigma _{j}$-group or 
a  $\sigma _{i}\cup \sigma _{j}$-group for
some   $j\ne i$. }

(iii) {\sl $G\simeq G^{*}/R$, where $G^{*}\leq A_{1}\times \cdots \times A_{t}$ for some 
$\sigma$-biprimary $\sigma$-soluble  groups  $A_{1}, \ldots ,  A_{t}$. }

Note that Satz 3.1  in \cite[VI]{hupp} can be obtained as 
   a special case of Theorem B,  when $\sigma =\{\{2\}, \{3 \}, \ldots \}$.

Finally, we prove the following 

{\bf Theorem E.}  
   {\sl  If  $G$ possesses  a complete Hall $\sigma$-set  $\cal  H$  with  
  $|G:N_{G}(H)|$   is $\sigma$-primary for all $H\in \cal H$, then
 $G$ is generalized  $\sigma$-soluble.}

{\bf Corollary 1.11} (See Zhang \cite {zhang} or   Guo \cite[Theorem 3]{Gu}). 
   {\sl   If   for every   Sylow subgroup $P$ of $G$
  the number $|G:N_{G}(P)|$   is
  a prime power, then $G$ is soluble.}

 \section{Preliminaries}

 We use ${\mathfrak{S}}_{g\sigma}$    to denote the class  of all 
$\sigma$-soluble    groups.

The direct calculations show that the following lemma is true 

{\bf Lemma 2.1.} (i) {\sl   The class  ${\mathfrak{S}}_{g\sigma}$ is
 closed under taking  
products of normal subgroups, homomorphic images and  subgroups.
 Moreover, any extension of the  generalized 
 $\sigma$-soluble group by a  generalized 
 $\sigma$-soluble group is  generalized 
 a $\sigma$-soluble group as well.    }

(ii) {\sl If $G/R, G/N \in {\mathfrak{S}}_{g\sigma}$, then $G/R\cap N
 \in {\mathfrak{S}}_{g\sigma}$.}

(iii)  {\sl ${\mathfrak{S}}_{g\sigma} \subseteq {\mathfrak{S}}_{g\sigma ^{*}}$ for any
 partition 
${\sigma ^{*}}=\{\sigma^{*}_{j}\ |\ j\in J\}$ of   $\Bbb{P}$ such that 
$J\subseteq I$ and $\sigma _{j}\subseteq \sigma ^{*}_{j}$ for all $j\in J$.}

Recall that $G$ is said to be: a $D_{\pi}$-group if $G$ possesses a Hall 
$\pi$-subgroup $E$ and every  $\pi$-subgroup of $G$ is contained in some 
conjugate of $E$;  a \emph{$\sigma$-full group
 of Sylow type} \cite{1} if every subgroup of $G$ is a $D_{\sigma _{i}}$-group for every
 $\sigma _{i}\in \sigma$.
  
In view of Theorem B in \cite{2}, the following fact is true.  

{\bf Lemma 2.1.} {\sl If $G$ is $\sigma$-soluble, then $G$ is a $\sigma$-full group
 of Sylow type.  
}

Let $\Pi \subseteq \sigma$. 
A natural number $n$ is said to be a \emph{$\Pi$-number}
 if  $\sigma (n)\subseteq \Pi$.  A  subgroup $A$ of $G$ is said to be: a 
\emph{ $\Pi$-subgroup} of $G$ if $\sigma (G)\subseteq \Pi$;
   a \emph{Hall $\Pi$-subgroup} of $G$ \cite{1} if    $|A|$ is  a $\Pi$-number 
 and $|G:A|$ is a $\Pi'$-number
   
We use  $\frak{H}_{\sigma}$ to denote the class of all 
groups $G$ such that  $G$ has a   complete Hall 
  $\sigma$-set ${\cal H}=\{H_{1}, \ldots , H_{t} \}$ satisfying the condition  
$H_{i}^{x}H_{j}^{y}=H_{j}^{y}H_{i}^{x}$  for all $x, y\in G$ and all 
$i\ne j$.  Note that in view of Lemma 2.1, each $\sigma$-biprimary 
$\sigma$-soluble  group  belongs to the class $\frak{H}_{\sigma}$.

 {\bf Lemma 2.2.} {\sl The class
$\frak{H}_{\sigma}$ is closed under taking homomorphic images and 
and direct products.  Moreover, if $G$ is a $\sigma$-full group
 of Sylow type  and $E\leq G\in \frak{H}_{\sigma}$, then $E\in \frak{H}_{\sigma}$.}

{\bf Proof.}   
 Let    $R$ be a normal subgroup 
of $G\in \frak{H}_{\sigma}$.  
By hypothesis,  $G$ has a   complete Hall 
  $\sigma$-set ${\cal H}=\{H_{1}, \ldots , H_{t} \}$ satisfying the condition  
$H_{i}^{x}H_{j}^{y}=H_{j}^{y}H_{i}^{x}$  for all $x, y\in G$ and all 
$i\ne j$.

Then ${\cal H}_{0}=\{H_{1}R/R, \ldots , H_{t}R/R \}$ is a         
 complete Hall   $\sigma$-set    of $G/R$ such that   
$$ (H_{i}R/R)^{xR}(H_{j}R/R)^{yR}=
      H_{i}^{x}H_{j}^{y}R/R=H_{j}^{y}H_{i}^{x}R/R$$$$=(H_{j}R/R)^{yR}(H_{i}R/R)^{xR}.$$
Thus  $G/R\in 
\frak{H}_{\sigma}$.  Therefore the class  $\frak{H}_{\sigma}$ is closed 
under taking homomorphic images.

Now   we show that the class
$\frak{H}_{\sigma}$ is closed under taking  direct products. It is enough 
to show that if 
$A, B \in \frak{H}_{\sigma}$, then $G=A\times B \in 
\frak{H}_{\sigma}$. Let ${\cal A}=\{A_{1}, \ldots , A_{n}\}$  be  a complete Hall   $\sigma$-set of $A$
 and  ${\cal B}=\{B_{1}, \ldots , B_{m}\}$ be  a complete Hall   $\sigma$-set of $B$ such that
 $A_{i}^{a_{1}}A_{j}^{a_{2}}=A_{j}^{a_{2}}A_{i}^{a_{1}}$  for all $a_{1}, a_{2}\in A$ and all 
$i\ne j$  and   $B_{i}^{b_{1}}B_{j}^{b_{2}}=B_{j}^{b_{2}}B_{i}^{b_{1}}$  for all
 $b_{1}, b_{2}\in B$ and all 
$i\ne j$.  We can assume  without loss of generality that  $1 \in {\cal 
A}\cap {\cal B} $ and that $\sigma (G)= \{\sigma _{1}, \ldots \sigma 
_{t}\}$. Therefore for every $i$ there are indices $a_{i}$ and $b_{i}$  such that 
$A_{a_{i}} \times B_{b_{i}}$  is a Hall  $\sigma _{i}$-subgroup of $G$.
Moreover, if $x=a_{1}b_{1}$ and $y=a_{2}b_{2} $, where $a_{1}, a_{2}\in A$ 
and   $b_{1}, b_{2}\in B$, then   $$(A_{a_{i}}
 \times B_{b_{i}})^{x}(A_{a_{j}} \times B_{b_{j}})^{y}=(A_{a_{i}}^{a_{1}}
 \times B_{b_{i}}^{b_{1}})(A_{a_{j}}^{a_{2}} \times B_{b_{j}}^{b_{2}})$$$$=
(A_{a_{j}}^{a_{2}} \times B_{b_{j}}^{b_{2}})(A_{a_{i}}^{a_{1}}
 \times B_{b_{i}}^{b_{1}})=(A_{a_{j}} \times B_{b_{j}})^{y}(A_{a_{i}}
 \times B_{b_{i}})^{x}.$$   Hence 
 $A_{a_{1}} \times B_{b_{1}}, \ldots , A_{a_{t}} \times B_{b_{t}}$ is a
 $\sigma$-basis of $G$. 
 Thus  $G \in \frak{H}_{\sigma}$.

Finally, assume that $G$ is a $\sigma$-full group
 of Sylow type. Then for every  complete Hall 
$\sigma$-set  ${\cal E}=\{E_{1}, \ldots , E_{r} \}$  of $E$ there is a    complete Hall
$\sigma$-set  ${\cal H}=\{H_{1}, \ldots , H_{t} \}$  of $G$  such that 
$E_{i}=H_{i}\cap E$ for all  $i=1, \ldots t$. We can assume without loss 
of generality that $E_{i}$ is a $\sigma _{i}$-group.

Now let  $\Pi =\{\sigma _{i}, \sigma _{j}\}$. 
Then for $x, y\in E$ we have 
  $\langle E_{i}^{x}, 
E_{j}^{y} \rangle \leq E\cap H_{i}^{x}H_{j}^{y}$, where   $E\cap 
H_{i}^{x}H_{j}^{y}$ is a $Pi$-subgroup of $E$ and so $|E\cap 
H_{i}^{x}H_{j}^{y}|\leq |E_{i}^{x}|| 
E_{j}^{y}|$.  Hence   $\langle E_{i}^{x}, 
E_{j}^{y} \rangle=E_{i}^{x}E_{j}^{y}$.  
Thus $E \in \frak{H}_{\sigma}$.

The lemma is proved.

The next lemma is evident.

 {\bf Lemma 2.3.} {\sl If the chief factors $H/K$ and $T/L$ of $G$ are $G$-isomorphic, then 
$C_{G}(H/K)=C_{G}(T/L)$.}

  We use  $\frak{X}_{\sigma}$ to denote  
the class of all  generalized $\sigma$-soluble groups  $G$ such that 
for  the automorphism group $G/C_{G}(H/K)$,  induced by $G$ on any  its 
chief factor $H/K$,  we have $|\sigma (G/C_{G}(H/K))|\leq 2$ and also 
 $\sigma(H/K)\subseteq  \sigma (G/C_{G}(H/K))$ in the case $|\sigma (G/C_{G}(H/K))|= 2$.

  {\bf Lemma 2.4.} {\sl The class
$\frak{X}_{\sigma}$ is closed under taking homomorphic images,  direct products
and subgroups.}

{\bf Proof.} Let    $R$ be a normal subgroup 
of $G\in \frak{X}_{\sigma}$.  Then $G/R$ is generalized $\sigma$-soluble  and
 for any chief factor $(H/R)/(K/R)$ of 
$G/R$ we have $C_{G/R}((H/R)/(K/R))=C_{G}(H/K)/R$, where $H/K$ is a chief 
factor of $G$. Hence   $G/R\in \frak{X}_{\sigma}$, so 
 the class  $\frak{X}_{\sigma}$  is closed 
under taking  homomorphic images.

Now let $G=A_{1}\times A_{2}$, where   $A_{1}, A_{2}\in \frak{X}_{\sigma}$.  The
 Jordan-H\"{o}lder 
theorem for groups with operators \cite[A, 3.2]{DH} implies that for every chief factor $H/K$
 of $G$, there is $i$ such that some chief factor $T/L$ of $G$ below $A_{i}$
 is $G$-isomorphic to $H/K$. Moreover,   $T/L $ is a chief factor of 
$A_{i}$ and $G/C_{G}(T/L)\simeq  A_{i}/C_{A_{i}}(T/L)$.  Hence, by Lemma 2.3, 
   $G\in \frak{X}_{\sigma}$, so   the class  $\frak{X}_{\sigma}$  is closed 
under taking   direct products.

 Finally, we show  that if  $G\in \frak{X}_{\sigma}$  and 
$E\leq G$, then $E\in \frak{X}_{\sigma}$. 

Let   $1=G_{0} < G_{1} < \cdots < G_{t-1} < G_{t}=G$ be a chief series of $G$. 
Let $H/K$ be any chief factor
of $E$ such that $G_{l-1}\cap E \leq K < H \leq G_{l}\cap E$ for some $l$.

Since  $$C_{G}(G_{l}/G_{l-1})\cap
E\leq C_{E}(G_{l}\cap E/G_{l-1}\cap E)\leq C_{E}(H/K),$$
$|\sigma (E/C_{E}(H/K))|\leq 2$. Moreover, if  $|\sigma 
(E/C_{E}(H/K))|= 2$, then  $|\sigma 
(G/C_{G}(G_{l}/G_{l-1}))|= 2$,   $\sigma 
(G/C_{G}(G_{l}/G_{l-1}))=\sigma (E/C_{E}(H/K)))$
 and $\sigma(G_{l}/G_{l-1})\subseteq 
 \sigma (G/C_{G}(G_{l}/G_{l-1}))$.  Hence from the isomorphism 
 $E\cap G_{l}/E\cap G_{l-1}\simeq (E\cap G_{l})G_{l-1}/G_{l-1}$ we get 
that $\sigma(H/K)\subseteq 
 \sigma (E/C_{E}(H/K))$.     Now
applying the Jordan-H\"older theorem for  
 groups with operators \cite[A, 3.2]{DH}  we get that $E\in 
\frak{X}_{\sigma}$.  Therefore the class  $\frak{X}_{\sigma}$  is closed 
under taking  subgroups.

The lemma is proved.

A \emph{class} of groups  is a collection $\mathfrak{X}$ of groups  with 
the property that if  $G\in \mathfrak{X}$  and if $H\simeq G$, then  
$H\in \mathfrak{X}$.   The symbol $(\mathfrak{Y})$  \cite[p. 264]{DH} denotes the smallest 
class of groups containing   $\mathfrak{Y}$. 

For a class $\mathfrak{X}$ of groups we define, following \cite[p. 264]{DH}: 

$\text{S}(\mathfrak{X})=(G: G\leq H$ for some $H\in \mathfrak{X});$

$\text{Q}(\mathfrak{X})=(G: G$ is an  epimorphic  image of some $H\in \mathfrak{X});$

$\text{D}_{0}(\mathfrak{X})=(G: G=H_{1}\times \cdots \times H_{r}$ for some 
 $H_{1},  \ldots     , H_{r}\in \mathfrak{X});$

$\text{R}_{0}(\mathfrak{X})=(G: G$ has normal subgroups $N_{1},  \ldots , 
N_{r}$ with   all  $G/N_{i}\in \mathfrak{X})$
 and $N_{1}\cap \cdots \cap N_{r}=1 )$.
 
We say, following \cite{shem}, that  a class $\mathfrak{X}$ of groups is a 
\emph{semiformation} if 
$$\text{S}(\mathfrak{X})=\mathfrak{X}=\text{Q}(\mathfrak{X}).$$

The following fact is well-known (see, for example, \cite[p. 57]{malc}).

{\bf Lemma  2.5.} {\sl  If $A\simeq B\leq G$, then for some $G^{*}\simeq G$ we 
have $A\leq  G^{*}$.}

{\bf Lemma 2.6.} {\sl  Let  $\mathfrak{X}$ be  a class  of groups.}

(1) {\sl If $\mathfrak{X}$ is a semiformation, then}
 $\text{S}(\text{D}_{0}(\mathfrak{X}))\text{=}R_{0}(\mathfrak{X}).$

(2) $\text{Q}(\text{R}_{0}(\mathfrak{X}))\subseteq \text{R}_{0}(\text{Q}(\mathfrak{X}))$.

(3)   $\text{R}_{0}(\text{R}_{0}(\mathfrak{X})) = \text{R}_{0}(\mathfrak{X})$.

{\bf Proof.}  (1) Let 
$\mathfrak{M} =\text{S}(\text{D}_{0}(\mathfrak{X}))$ and $\mathfrak{H} = 
\text{R}_{0}(\mathfrak{X})$. First   suppose that $G\in  \mathfrak{M}$. 
Then,  in view of Lemma 2.5,  $G\leq  H_{1}\times \cdots \times H_{r}$ for some   
 $H_{1},  \ldots     , H_{r}\in \mathfrak{X}$.  We show that  $G\in  
\mathfrak{H}$. If $r=1$, it is clear. Now assume that $r > 1$ and let
 $N_{i}=H_{1}\times
  \cdots H_{i-1}H_{i+1}\times \cdots  \times H_{r}$ for all $i=1, \ldots , r$. 
 Then $$G/G\cap N_{i}\simeq GN_{i}/N_{i} \leq   (H_{1}\times \cdots \times 
H_{r})/N_{i}\simeq H_{i} \in  \mathfrak{X}.$$  It is clear also that $G\cap N_{1}, \ldots G\cap N_{r}$  
are normal subgroups of $G$ with   $(G\cap N_{1})\cap  \cdots \cap (G\cap N_{r})=1$. Hence   
$G\in \mathfrak{H}$.    Therefore   $\mathfrak{M} \subseteq  \mathfrak{H}$.

Now  suppose that $G\in  \mathfrak{H}$, that is, $G$ has normal subgroups $N_{1},  \ldots , 
N_{r}$ with   all  $G/N_{i}\in \mathfrak{X}$
 and $N_{1}\cap \cdots \cap N_{r}=1$.  Then, by Lemma 4.17 in \cite[A]{DH}, 
 $G$ is isomorphic with  a subgroup of
  $ (G/N_{1})\times \cdots \times (G/N_{r}).$   Hence 
$G\in \mathfrak{M}$, so
 $\mathfrak{H} \subseteq  \mathfrak{M}$.
Therefore we have (1).

(2), (3) See respectively Lemmas   1.18(b) and Lemma 1.6  in \cite[II]{DH}.

The lemma is proved.

{\bf Lemma 2.7.} {\sl Let   $\mathfrak{X}$  be a  semiformation of groups. Suppose that 
$\frak{F}$ is the class of groups $A$  which can be represented in the form
 $A\simeq A^{*}/R$, where $A^{*}\leq A_{1}\times \cdots \times A_{t}$ for some 
 $A_{1}, \ldots ,  A_{t}\in \frak{X}$. If  $G/R, G/N\in\frak{F}$, then  
$G/(R \cap N)\in\frak{F}$.
    }

{\bf Proof. }     We can assume without loss of generality that $R\cap N=1$.
 Then, by Lemma 2.5, 
 $$G\in \text{R}_{0}(\text{Q}(\text{S}(\text{D}_{0}(\mathfrak{X}))))=
\text{R}_{0}(\text{Q}(\text{R}_{0}(\mathfrak{X})))\subseteq
 \text{Q}(\text{R}_{0}(\text{R}_{0}(\mathfrak{X})))=\text{Q}(\text{R}_{0}(\mathfrak{X})
=\text{Q}(\text{D}_{0}(\mathfrak{X})))\subseteq \frak{F}.$$ 

  The lemma is proved.

{\bf Lemma 2.8.} {\sl Let   $G=RH_{1} \cdots H_{n}$, where
 $[H_{i}^{G}, H_{j}^{G}]=1$ for all $i\ne j$ and $R$ is normal in $G$.
 Then $G\simeq G^{*}$, where  $G^{*}\leq (RH_{1})\times \cdots \times (RH_{n})$.   }

{\bf Proof.}  See pages 671--672  in \cite{hupp}.

{\bf Lemma 2.9} (See  \cite[2.2.8]{15}). {\sl If $\frak{F}$ is a non-empty 
formation and 
$N$, $R$ be  subgroups of $G$, where $N$ is normal in $G$.}

(i) {\sl 
 $(G/N)^{\frak{F}}=G^{\frak{N}}N/N.$  }

(ii) {\sl If $G=RN$, then $G^{\frak{N}}N=R^{\frak{N}}N$}.

{\bf Lemma 2.10.} {\sl  Suppose that $G$ has a  $\sigma _{i}$-subgroup 
$A\ne 1$ and a  $\sigma _{j}$-subgroup $B\ne 1$   such that 
$AB^{x}=B^{x}A$ for all $x\in G$. If $O_{\sigma _{i}\cup \sigma _{j}}(G)= 1$, then
 $[A^{G}, B^{G}]= 1$.  
   }

{\bf Proof.} By hypothesis,
$A(B^{x})^{y}=(B^{x})^{y}A$ for all $x, y\in G$ and    $$D=\langle
A^{B^{x}} \rangle \cap \langle (B^{x})^{A}\rangle\leq  AB^{x},$$ where $D$ is
subnormal in $G$ by \cite[1.1.9(2)]{prod}. Then $D\leq O_{\sigma _{i}\cup \sigma 
_{j}}(G)= 1$, so 
Then $ [A, B^{x}] \leq  [\langle A^{B^{x}}\rangle , \langle B^{x})^{A}\rangle ]\leq  D = 1$.
Therefore  $[A^{x}, B^{y}]=1$ for all $x, y\in G$ and so $[A^{G}, B^{G}]= 1$.
 The lemma is proved.

\section{Proofs of the results}

 {\bf Proof of Theorem B.}   
 Suppose  that this theorem  is false and let $G$ be a
counterexample with $|G|$ minimal.  We can assume without loss of        
generality   that $i=1$, $j=2$ and $k=3$.

(1) {\sl If $R$ is a minimal normal subgroup of $G$, then  $G/R$ is generalized
 $\sigma$-soluble.
 Hence $R$ is the unique minimal normal subgroup of $G$ and $R$ is not generalized
 $\sigma$-soluble.
 Thus  $C_{G}(R)=1$. }

First note that if $A_{i}\leq R$, then for $j\ne i$ we have $G= A_{i}A_{j} 
=RA_{j}$, so $G/R\simeq A_{j}/A_{j}\cap R$  is generalized
 $\sigma$-soluble.

Now suppose that $A_{i}\nleq R$ for all $i=1, 2, 3$. 
 For any group $A$ and any $l\in I$ we have
 $O^{\sigma _{l}}(A)=A^{\frak{G}_{\sigma _{l}}}$, where   
$\frak{G}_{\sigma _{l}}$ is the class of all  $\sigma _{l}$-groups. 
Therefore  $O^{\sigma _{l}}(A)R/R=O^{\sigma _{l}}(AR/R)$   by Lemma 2.10(ii). 
Hence 

  $$N_{G}(O^{\sigma _{l}}(A))R/R\leq 
N_{G}(O^{\sigma _{l}}(A)R)/R  =  N_{G/R}(O^{\sigma _{l}}(AR/R).$$ 
 Hence   the  three  indices  
$|(G/R):N_{G/R}((O^{\sigma _{1}}(A_{1}))|$,  $
|(G/R):N_{G/R}(O^{\sigma _{2}}(A_{2}))|,$ and 
  $ |(G/R):N_{G/R}(O^{\sigma _{3}}(A_{3}))|$ are pairwise  $\sigma$-coprime.
 On the other hand, clearly, $A_{1}R/R\simeq A_{1}/A_{1}\cap R$ is  $\sigma$-soluble
 and
  $A_{i}R/R\simeq A_{i}/A_{i}\cap R$ is generalized $\sigma$-soluble
 for all $i=2, 3$.
 Therefore  the hypothesis holds
 for $G/R$, so $G/R$ is generalized $\sigma$-soluble
 by the choice of $G$. Hence $R$ is  not generalized  $\sigma$-soluble
 and 
 $R$ is the unique minimal normal subgroup
 of $G$ by Lemma 2.1(ii). Therefore, since $C_{G}(R)$ is normal in $G$, $C_{G}(R)=1$.

(2) {\sl At least two of
 the subgroups $A_{1}, A_{2},  A_{3}$ are not $\sigma$-primary.}

Indeed, if $A_{i}$ and   $A_{j}$  are  $\sigma$-primary for some $i\ne j$,
 then $G=A_{i}A_{j}$ is generalized $\sigma$-soluble, which contradicts the choice of $G$.

(3) {\sl  If $A_{l} $ is  $\sigma$-soluble and $L$ is a minimal normal 
subgroup of $A_{l}$, then for some $t\ne l$  we have  $O^{\sigma 
_{t}}(A_{t})=1$, so  $A_{t}$ is $\sigma$-primary.}

We can assume without loss of generality  that  $l=2$. 
Let $L$ be   a minimal  normal subgroup  of $A_{l}=A_{2}$.  Then  $L$
 is a $\sigma _{i}$-group for
 some prime $i$ since  $A_{2}$ is $\sigma$-soluble by hypothesis.
  On the other                                                   
 hand, again by hypothesis,    
$|G:N_{G}(O^{\sigma _{1}}(A_{1}))|$ and 
 $|G:N_{G}(O^{\sigma _{3}}(A_{3}))|$ are   $\sigma$-coprime, so  
 at least one of the numbers,
 $|G:N_{G}(O^{\sigma _{3}}(A_{3}))|$ say,   is a  $\sigma _{i}'$-number. 
But $G=A_{2}A_{3}=A_{2}N_{G}(O^{\sigma _{3}}(A_{3}))$, so 
 $|A_{2}:A_{2}\cap N_{G}(O^{\sigma _{3}}(A_{3}))|$    is a  $\sigma _{i}'$-number. 
Hence   $$L\leq A_{2}\cap N_{G}(O^{\sigma _{3}}(A_{3}))\leq N_{G}(O^{\sigma _{3}}(A_{3})).$$
 Therefore 
$$R\leq L^{G}=L^{A_{2}N_{G}(O^{\sigma _{3}}(A_{3}))}=
L^{N_{G}(O^{\sigma _{3}}(A_{3}))}\leq N_{G}(O^{\sigma _{3}}(A_{3}))$$   by 
Claim (1). 

Suppose that  $O^{\sigma _{3}}(A_{3})\ne 1$  and let  $H=O^{\sigma _{3}}(A_{3}) \cap R$. 
If $H\ne 1$, then $H$ is a non-identity normal  generalized $\sigma$-soluble 
subgroup of $R$ since $A_{3}$ is  generalized $\sigma$-soluble  by hypothesis.  Hence $R$ is generalized  $\sigma$-soluble, contrary to Claim (1).
Thus   $H=1$, so  $O^{\sigma _{3}}(A_{3})\leq 
C_{G}(R)\leq R$, so $O^{\sigma _{3}}(A_{3}) =1$ by Claim (1).  This 
contradiction completes the proof of (3).

{\sl Final contradiction.} Since $A_{1}$ is   $\sigma$-soluble by 
hypothesis, Claim (3) implies that one of the subgroups $A_{2}$ or $A_{3}$, $A_{2}$ say, is
 $\sigma$-primary. Then  $A_{2}$ is   $\sigma$-soluble   and so, again by Claim (3), one  
 of the subgroups $A_{1}$ or $A_{3}$ is also $\sigma$-primary.    But this 
is impossible by Claim (2). 
  This contradiction completes the proof of 
the result.

{\bf Proof of Theorem A.}  
  We 
can assume without loss of generality that $H$ is a $\sigma _{1}$-group, $2\in \sigma _{1}$   
and   $ {\cal H}=\{H_{1}, \ldots , H_{t}\}$, where $H_{i}$ is a $\sigma _{i}$-group for all
 $i=1, \ldots, t$.  Then $ t > 2$ since otherwise $G$ and $H^{G}$ are 
generalized  $\sigma$-soluble.

 (i) This assertion, in fact, is a corollary of Theorem D. Indeed,  
 let 
 $E_{i}=H_{1}\cdots H_{i-1}H_{i+1} \cdots H_{t}$  
  for all $i= 1, \ldots , t$. Then,  by 
induction, $E_{1}, \ldots , E_{t}$ are  generalized  $\sigma$-soluble and, by the
 Feit-Thompson theorem, $E_{1}=H_{2}\cdots  H_{t}$ is soluble. Since $t > 
2$, and evidently,  the indeces $|G:E_{1}|$, $|G:E_{2}|$, 
 $|G:E_{3}|$ are pairwise  $\sigma$-coprime, $G$ is   generalized 
 $\sigma$-soluble by  Corollary  1.6.

 (ii) Assume that this assertion  is falls and let $G$ be a
 counterexample of minimal order.    By hypothesis, $HA^{x}=A^{x}H$ 
for all $x\in G$ and all  $A  \in  \cal H$ such that $(|H|, |A|)=1$.

(1) {\sl For some $i > 1$ and $x\in G$  we have $H_{i}^{x}\nleq N_{G}(H)$.}

Indeed, suppose that for all  $i > 1$ and all $x\in G$  we have $H_{i}^{x}\leq 
N_{G}(H)$.   Then $$E=(H_{2})^{G} \cdots (H_{t})^{G}\leq N_{G}(H)$$
and so $H^{G}=H^{H_{1}E}=H^{H_{1}}\leq H_{1}$ is generalized 
$\sigma$-soluble, a contradiction. Hence we have (1).

(2) {\sl $H^{G}R/R$ is generalized  $\sigma$-soluble. }

It is clear that $ {\cal H}_{0}=\{H_{1}R/R, \ldots , H_{t}R/R\}$  is a 
complete Hall $\sigma$-set of $G/R$ and $HR/R\leq   H_{1}R/R$. Moreover, 
$$(HR/R)(H_{i}R/R)^{xR}=HH_{i}^{x}R/R=H_{i}^{x}HR/R=(H_{i}R/R)^{xR}(HR/R) $$ for all 
  $ i > 1$ and all $xR\in G/R$. Therefore  $HR/R$ 
 is ${\sigma}$-semipermutable
 in $G/R$                                     
 with respect to ${\cal H}_{0}$, so $(HR/R)^{G}=H^{G}R/R$ is  generalized  $\sigma$-soluble 
by  the choice of $G$.

(3) {\sl Every minimal normal subgroup of $G$ is not generalized  $\sigma$-soluble. }
                                         
Indeed, if $R$ is  generalized  $\sigma$-soluble, then  from the isomorphism $H^{G}R/R\simeq 
H^{G}/H^{G}\cap R$ and Claim (2) we get that $H^{G}$ is
  generalized  $\sigma$-soluble, contrary to the choice of $G$.   

 {\sl Final contradiction for (ii).}  Let $i > 1$. By hypothesis,
$HH_{i}^{x}=H_{i}^{x}H$  for all $x\in G$. Om the other hand,  $O_{\sigma 
_{1}\cup \sigma _{i}}(G)=1$ by Claim (3). Therefore $[H^{G}, H_{i}^{G}]= 1$ by Lemma 2.11, contrary to Claim (1).  Hence 
Assertion (ii) holds.

The theorem is proved.

In fact, Theorem C is a corollary of the following

{\bf Proposition 3.1.} {\sl If $G$ is $\sigma$-full, then
 $ \frak{H}_{\sigma}=\frak{X}_{\sigma} $.}

{\bf  Proof.}  First we show $ \frak{H}_{\sigma}\subseteq \frak{X}_{\sigma} $.
  Assume that this is false and let $G$ be a group  of
 minimal orderin $ \frak{H}_{\sigma}\setminus  \frak{X}_{\sigma} $.
  Then $|\sigma (G)| > 2$.  By hypothesis,  $G$ has a   complete Hall 
  $\sigma$-set ${\cal H}=\{H_{1}, \ldots , H_{t} \}$ satisfying the condition  
$H_{i}^{x}H_{j}^{y}=H_{j}^{y}H_{i}^{x}$  for all $x, y\in G$ and all 
$i\ne j$.    
 We can assume without loss of 
generality that $H_{i}$    is a $\sigma _{i}$-group for all $i=1, \ldots , t$. 
Let $R$ be a minimal normal subgroup of $G$

(1) {\sl $G/R\in \frak{X}_{\sigma}$.  Hence $R$ is the only minimal normal subgroup of $G$.  }

First note that since  $G\in \frak{H}_{\sigma}$, we have  $G/R\in \frak{H}_{\sigma}$ by Lemma 
2.2 and so
  $G/R\in \frak{X}_{\sigma}$ by the choice of $G$. Now assume that $G$ has a minimal normal 
subgroup $N\ne R$.  Then $G/N\in \frak{X}_{\sigma}$, so from
 the $G$-isomorphism $RN/N\simeq R$, Lemma 2.3  
and the Jordan-H\"older theorem for  
 groups with operators \cite[A, 3.2]{DH}  we get  that  $G\in \frak{X}_{\sigma}$, a 
contradiction.  Hence  we have (1).

(2) {\sl $G$ is generalized $\sigma$-soluble. Hence $R$ is
 a $R$ is a $\sigma _{i}\cup \sigma _{j}$-group for some $i, j$}
 (This follows from Theorem C(i)).

(3) {\sl $O_{\sigma _{l}\cup \sigma _{k}}(G)= 1$ for all $k, l$ such that
 $\sigma (R) \not \subseteq O_{\sigma _{i}, \sigma _{j}}$. } (This
 follows from Claims (1) and (2)).

(4)   $ \frak{H}_{\sigma}\subseteq \frak{X}_{\sigma} $.

 Let 
$\sigma _{k}\in \sigma (G)$  such that $i\ne k\ne j$. 
First assume that $\sigma (R)=\{\sigma _{i},  \sigma _{j}\}$. Then 
$R=(R\cap H_{i})(R\cap H_{j})$.  Claim (1) implies that      $O_{\sigma 
_{i}\cup \sigma _{k}}(G)= 1$ and  $O_{\sigma 
_{j}\cup \sigma _{k}}(G)= 1$. Then, by Lemma 2.11, $H_{k}^{G}\leq C_{G}(H_{i})$ and  
$H_{k}^{G}\leq C_{G}(H_{j})$.    Therefore  $H_{k}^{G}\leq C_{G}(R)$. 
Therefore $G/C_{G}(R)$ is a $\sigma _{i}\cup \sigma _{j}$-group and so 
$G\in    \frak{X}_{\sigma}$ since  $G/R\in \frak{X}_{\sigma}$ by Claim (1), a contradiction.
Hence  we can assume that   $\sigma (R)=\{\sigma _{i}\}$. 

Let $l\ne k$ and   $l\ne i\ne k$. Then  $O_{\sigma _{l}\cup \sigma 
_{k}}(G)= 1$ by Claim (3). Hence $[H_{l}^{G}, H_{k}^{G}]=1$ by Lemma 2.11.
 Therefore from Claim (1) we get that $H_{l}\leq  C_{G}(R)$ for all $l\ne 
i$ and  so $G/C_{G}(R)$ is a $\sigma _{i}$-group.  Therefore  we again get 
that $G\in    \frak{X}_{\sigma}$.   
This contradiction completes the proof of the inclusion
 $ \frak{H}_{\sigma}\subseteq \frak{X}_{\sigma} $.

Now we show that if $G\in \frak{X}_{\sigma}$, then
   $VW=WV$ for every Hall $\sigma _{i}$-subgroup $V$, 
 every Hall $\sigma _{j}$-subgroup $W$ of $G$ and all $i\ne j$. 
Assume  that this is false and let $G$ be a counterexample of minimal 
order.   Then $|\sigma (G)| > 2$. 
  
(a) {\sl   $G=RVW$.  Hence $G/R$ is a $\sigma _{i}\cup \sigma _{j}$-group}

Lemma 2.4 implies that the hypothesis holds for $G/R$, so the choice
 of $G$ implies that   $$(VR/R)(WR/R)=(WR/R)(VR/R)=VWR/R.$$ Hence   $RVW$ is a subgroup
of $G$.  Suppose that  $RVW < G$.  For each $k\ne i, j$, $H_{k}\cap R$ is
 a Hall $\sigma _{k}$ subgroup of $RVW$.  
Therefore the hypothesis holds for $RVW$ by Lemma 2.4 and so 
 $VW=WV$, a contradiction. Thus we have (a).

(b) {\sl $R$ is the unique minimal normal subgroup of $G$. Hence $R$ is a $p$-group for some
 $p\in \sigma _{i}$. }

Indeed, suppose that $G$ has a  minimal normal subgroup  $N\ne R$. Then
 $G/N$ is a $\sigma _{i}\cup \sigma _{j}$-group by Claim (1), 
so $G\simeq G/R\cap N$  is a $\sigma _{i}\cup \sigma _{j}$-group. Hence 
$|\sigma (G)|=2$, a contradiction. Therefore $R$ is the unique minimal normal 
subgroup of $G$. Suppose that $R$ is non-abelian. Then $C_{G}(R)=1$,  
since $C_{G}(R)$  is normal in $G$. Hence $G\simeq G/C_{G}(R)$ is $\sigma 
$-biprimary, a contradiction. Thus we have (b).

(c) {\sl     $R$ is   a Sylow $p$-subgroup of $G$.}

 Assume that $R$ is not a Sylow $p$-subgroup of $G$. Since $G=RVW$, by 
Claim (5), but $VW\ne WV$, it follows that $R\nleq VW$. Hence for a Sylow 
$p$-subgroup $P$ of $G$ we have $P\cap V=1=P\cap W$. Hence $P\nleq RVW=G$, 
a 
contradiction. Thus  we have (c).

(d)   $ \frak{X}_{\sigma}\subseteq \frak{H}_{\sigma} $.

In view of Claims  (b) and (c),  there is a
maximal subgroup $M$ of $G$ such that $G=R M$ and $M_{G}=1$.
Hence, for some $k$, we have  $R=C_{G}(R)=O_{p}=H_{k}$ by \cite[A, 15.6]{DH}.  
Then $\sigma _{k} \not \in \sigma (G/C_{G}(R))$, so $G/C_{G}(R))=G/R$ is a 
$\sigma _{l}$-group  for some $l$ since  $G\in {\frak{G}}_{\sigma}$ by 
hypothesis. But then $|\sigma (G)| \leq   2$.
 This contradiction completes the proof of the inclusion 
  $ \frak{X}_{\sigma}\subseteq \frak{H}_{\sigma} $.

From Claims (c) and (d) we get that   $ \frak{X}_{\sigma} = 
\frak{H}_{\sigma} $.   
The theorem is proved.

{\bf Proof of Theorem D.}  The implications   (i) $\Rightarrow$ (ii) and (ii) $\Rightarrow$ 
(i) directly follow from Theorem A.

Let   $\frak{X}$ 
 be the    class of  all $\sigma$-biprimary $\sigma$-soluble  groups.
 Then from Lemma 
2.2  it follows that if  $A$ is a group satisfying   $A\simeq A^{*}/R$ for  
 some $A^{*}\leq A_{1}\times \cdots \times A_{t}$ and  $A_{1}, \ldots ,  A_{t}\in \frak{X}$,
 then $A$  satisfies also Condition (i). 
Thus    (iii) $\Rightarrow$  (i). 

Now we show that  (i) $\Rightarrow$  (iii).  Assume that this is false and 
let  $G$ be a group of minimal order among the groups which satisfy 
Condition (i) but do  not satisfy Condition (iii). Since $G$ is $\sigma$-soluble, it has 
a complete Hall $\sigma$-set    $ {\cal H}=\{H_{1}, \ldots , H_{t}\}$ and, by hypothesis, 
$H_{i}^{x}H_{j}^{y}=H_{j}^{y}H_{i}^{x}$  for all $x, y\in G$ and all 
$i\ne j$.  We can assume without loss of generality that 
 $H_{i}$ is a $\sigma _{i}$-group for all
 $i=1, \ldots, t$.  Then $ t > 2$ since otherwise Condition (iii) holds 
for $G$.

 Let $L$ be a minimal normal 
subgroup of $G$.
 Since $G$ is $\sigma$-soluble, $L$ is 
$\sigma$-primary, $L$ is a $\sigma _{1}$-group say.  Let  
$\frak{F}$ be the class of groups $A$  which can be represented in the form
 $A\simeq A^{*}/R$, where $A^{*}\leq A_{1}\times \cdots \times A_{t}$ for some 
 $A_{1}, \ldots ,  A_{t}\in \frak{X}$.    
Lemma 2.2 implies that Condition (i) holds for $G/L$, so  $G/L \in 
\frak{F}$ by the choice of $G$. If $G$ has a minimal normal subgroup $N\ne 
L$, then also we have    $G/N \in \frak{F}$, so  $G\simeq G/1 =G/L\cap N \in \frak{F}$
 by Lemma 
2.7. Thus $L$ is the unique minimal normal subgroup of $G$.

Now let $j\ne i\ne 1\ne j$, and let $A=H_{i}$ and $B=H_{j}$. 
 Then
$AB^{x}=B^{x}A$ for all $x\in G$.  It is clear also that
 $O_{\sigma _{i}\cup \sigma _{j}}(G)=1$.  Then  $[A^{G}, B^{G}]= 1$ by Lemma 2.10.  
         Now using Lemma 2.8, we get that $G\in \frak{F}$. This
 contradiction completes the proof of the result.

{\bf Proof of Theorem E.}   Assume that this is false and let $G$ be a counterexample of
 minimal order.  Let    
${\cal H}=\{H_{1}, \ldots ,  H_{t} \}$.
Then $ t > 2$.  We can assume without loss of generality that $H_{i}$ is a 
non-identity $\sigma _{i}$-group   for all $i=1, 
\ldots , t$.     Let $N_{i}=N_{G}(H_{i})$ for all $i=1, 
\ldots , t$.

First we show that  $G$ is not simple. 
Assume that $G$ is a simple non-abelian group. Then  
$|G:N_{i}|\ne 1$ is  $\sigma$-primary for all
$i=1, \ldots , t$.

Let $|G:N_1|$
be a $\sigma_{i}$-number and  $I_{0}=\{1,  \ldots , i\} \setminus \{1, i\}$.
Then $I_{0} \ne \emptyset$ since $t > 2$. We can assume without loss of generality that
  $i=t$.
If for some
$k\in I_{0}$, $|G:N_{k}|$ is not a $\sigma_{t}$-number, then
$G=N_{1}N_{k}$ and  $|G:N_{1}\cap N_{k}|=|G:N_{1}||G:N_{k}|$  is a $\sigma 
_{k}'$-number, so $H_{k}\leq N_{1}\cap N_{k}$ since $N_{1}\cap N_{k}\leq N_{G}(H_{k})$. 
Therefore  $$H_{k}^{G}=  H_{k}^{N_{k}N_{1}}=H_{k}^{N_{1}}\leq N_{1},$$
 so $G$ is not simple  since $N_{1}\ne G$, a contradiction.   
 Therefore
$|G:N_{k}|$ is  a  $\sigma_{t}$-number for all  $k=1, \ldots , t-1$.

Let  $|G:N_{t}|$ be a  $\sigma_{k}$-number and $j\in \{1, \ldots , t\}\setminus
\{k, t \}$. Then  $G=N_{j}N_{t}$  and  $H_{j}\leq N_{j}\cap N_{t}$,
Hence  $H_{j}^{G}= H_{j}^{N_{j}N_{t}}=H_{j}^{N_{t}}\leq N_{t}$,
 so $G$ is not simple. 

Now let $R$ be any non-identity normal subgroup of $G$. Then 
       $\{H_{1}R/R, \ldots ,
H_{t}R/R \}$ is  a complete Hall $\sigma$-set  of
 $G/R$.  On the other hand, since $N_{i}R/R\leq 
N_{G/R}(H_{i}R/R)$ and $|G:N_{i}|$ is  $\sigma$-primary, 
$N_{G/R}(H_{i}R/R)$  is  $\sigma$-primary.  
Therefore the hypothesis holds on $G/R$ for any minimal normal
subgroup $R$ of $G$.   Hence   the choice of $G$ implies that $G/R$ is  
generalized 
$\sigma$-soluble, and consequently $R$ is not generalized 
$\sigma$-soluble.

  It is clear  that   $\{R_{1}, \ldots ,
R_{t}\}$ is  a complete Hall $\sigma$-set  of
$R$, where  $R_{i}=H_{i}\cap R$ for all $i=1, \ldots , t$.
Moreover, we have also $N_{i}\leq N_{G}(R_{i})$ and so from  $|N_{i}R:N_{i}|=|R:R\cap  N_{i}|$
 we get that  $|R:N_{R}(R_{i})|$   is $\sigma$-primary. 
 Therefore the hypothesis holds for $R$ and $R < G$ since $G$ is not 
simple.  Therefore $R$ is generalized $\sigma$-soluble   by the choice of 
$G$.  This contradiction completes the proof of the result.

\end{document}